\documentclass[12pt]{article}
\usepackage{amsmath}
\usepackage{amsfonts}
\usepackage{graphicx}
\usepackage{amssymb}
\usepackage{theorem}
\usepackage{tikz}
\usepackage{tikz-cd}

\usepackage[OT2,OT1]{fontenc}
\newcommand\cyr{%
\renewcommand\rmdefault{wncyr}%
\renewcommand\sfdefault{wncyss}%
\renewcommand\encodingdefault{OT2}%
\normalfont
\selectfont}
\DeclareTextFontCommand{\textcyr}{\cyr}

\newtheorem{theorem}{Theorem}
\newtheorem{definition}{Definition}
\newtheorem{corollary}{Corollary}

\newtheorem{lemma}{Lemma}

\newtheorem{remark}{Remark}

\newtheorem{construction}{Construction}
\usepackage[russian, english]{babel}
\usepackage{centernot}

\usepackage[unicode,colorlinks,plainpages=false]{hyperref}

\newcommand{\openbox}{$\begin{array}{c}
\hspace*{-0.55em}\sqcap \hspace*{-0.60em}\\[-0.4em] \hline
\multicolumn{1}{c}{\hspace*{-0.60em}}\\[-0.8em]
\end{array}$}

\begin{document}

\centerline{\bf On special Rees matrix semigroups over semigroups\footnote{Corresponding author: Attila Nagy,
email: nagyat@math.bme.hu
\\Researct is supported by National Research, Development and Innovation Office – NKFIH, 115288}}

\bigskip

\centerline{Attila Nagy and Csaba T\'oth}

\medskip

\centerline{Department of Algebra}
\centerline{Institute of Mathematics}
\centerline{Budapest University of Technology and Economics}
\centerline{1521 Budapest, P.O. Box 91, Hungary}
\centerline{nagyat@math.bme.hu, tcsaba94@gmail.com}
\bigskip

\begin{abstract}
In this paper we focus on Rees $I\times \Lambda$ matrix semigroups without zero over a semigroup $S$ with $\Lambda\times I$ sandwich matrix $P$, where $I$ is a singleton, $\Lambda$ is the factor semigroup of $S$ modulo the kernel $\theta_S$ of the right regular representation of $S$, and $P$ is a choice function on the collection of all $\theta _S$-classes of $S$. We describe the kernel of the right regular representation of this type of Rees matrix semigroups, and prove embedding theorems on them. Motivated by one of embedding theorems, we show how right commutative right cancellative semigroups can be constructed. We define the concept of a right regular sequence of semigroups, and show that every congruence on an arbitrary semigroup defines such a sequence.
\end{abstract}

\bigskip


\medskip


\bigskip

\section{Introduction}

The concept of the Rees matrix semigroup is one of the tools that can be effectively applied in the study of the structure of semigroups.
The Rees matrix construction was first used by D. Rees in \cite{Rees} to give a description of completely simple and completely $0$-simple semigroups in terms of their maximal subgroups.
Groups in the construction can be replaced with semigroups, see, for example, \cite{Ayik}, \cite{Descalco}, \cite{Kambites}, \cite{MacAlister}, \cite{Petrich}. If $S$ is a semigroup, $I$ and $\Lambda$ are non empty sets, and $P$ is a $\Lambda \times I$ matrix with entries $P(\lambda , i)$, then the set ${\cal M}(S; I, \Lambda ;P)$ of all triples $(i, s, \lambda )\in I\times S\times \Lambda$ is a semigroup under the multiplication $(i, s, \lambda )(j, t, \mu )=(i, sP(\lambda , j)t, \mu )$. This semigroup is called a \emph{Rees $I\times \Lambda$ matrix semigroup without zero over the semigroup $S$ with $\Lambda \times I$ sandwich matrix $P$}. If $I=\{ i\}$ is a singleton, then the entries of $P$ can be denoted by $P(\lambda)$, the element $(i, s, \lambda )$ of
${\cal M}(S; I,\Lambda ;P)$ can be considered in the form $(s, \lambda )$, and the operation on ${\cal M}(S; I, \Lambda ;P)$ is modified as follows:
$(s, \lambda )(t, \mu )=(sP(\lambda)t, \mu)$. In this case, we shall use the notation ${\cal M}(S;\Lambda ;P)$ instead of ${\cal M}(S; I, \Lambda ;P)$.
The concept of this type of Rees matrix semigroups is not new, it has already occurred in the book \cite{Clifford2:sg-20}.

Let $S$ be a semigroup and $\alpha$ be a congruences on $S$. Let $P$ be a choice function on the collection of all $\alpha$-classes of $S$, i.e., $P([s]_{\alpha})\in [s]_{\alpha}$ for every $\alpha$-class $[s]_{\alpha}$ of $S$.
Consider the Rees matrix semigroup ${\cal M}(S; S/\alpha; P)$.
Let $\beta$ be a congruence on $S$ with $\alpha \subseteq \beta$, and let
$P': S/\alpha \to S/\beta$ be the mapping defined as follows: $P'([s]_{\alpha})=[s]_{\beta}$ for every $[s]_{\alpha}\in S/\alpha$. Consider the Rees matrix semigroup ${\cal M}(S/\beta;S/\alpha;P')$.
Let
${\cal F}_{\alpha , \beta}$ be the mapping of ${\cal M}(S; S/\alpha; P)$ onto
${\cal M}(S/\beta; S/\alpha; P')$ defined as follows:
${\cal F}_{\alpha, \beta}:\ (a, [s]_{\alpha})\mapsto ([a]_{\beta}, [s]_{\alpha})$
for every $(a, [s]_{\alpha})\in {\cal M}(S; S/\alpha; P)$. In Section 2, we prove that ${\cal F}_{\alpha, \beta}$ is a surjective homomorphism.

\medskip

If $a$ is an arbitrary element of a semigroup $S$, then $\varrho _a:\ s\mapsto sa\ (s\in S)$ is a right translation (\cite{Clifford1:sg-1}) of $S$, and $\Phi _S :a\mapsto \varrho _a$ is a homomorphism of the semigroup $S$ into the semigroup of all right translations of $S$. The homomorphism $\Phi _S$ is called the \emph{right regular representation} of $S$. The right regular representation of semigroups plays an important role in the investigation of the structure of semigroups. Here we refer to only papers \cite{Chrislock1}, \cite{Chrislock2}, \cite{NagyKolibiar:sg-19}, \cite{Nagy:sg-8}, and \cite{Nagy:sg-20} in which the structure of semigroups $S$ are examined by the help of the semigroup $\Phi_S(S)$.
The right regular representation $\Phi _S$ is injective if and only if $S$ is \emph{left reductive}, i.e., for arbitrary elements $a, b\in S$, $sa=sb$ for all $s\in S$ implies $a=b$.
For an arbitrary semigroup $S$, let $\theta _S$ denote the \emph{kernel of $\Phi_S$}. It is clear that $(a,b)\in \theta _S$ if and only if $xa=xb$ for all $x\in S$. For every semigroup $S$, the factor semigroup $S/\theta _S$ will be denoted by $S/\theta$.

In this paper we focus on Rees matrix semigroups ${\mathcal M}(S;S/\theta;P)$, where $S$ is an arbitrary semigroup and $P$ is an arbitrary choice function on the collection of all $\theta_S$-classes of $S$.
In Section 2, we investigate the kernel of the right regular representation of ${\mathcal M}(S;S/\theta;P)$. In Section 3, we discuss the embedding of semigroups ${\mathcal M}(S;S/\theta ; P)$ in idempotent-free left simple semigroups.

Let $\varrho$ be a congruence on a semigroup $S$. Then $\varrho ^{*}=\{ (a, b)\in S\times S: (\forall s\in S) (sa, sb)\in \varrho\}$ is also a congruence on $S$ which will be called the \emph{right colon congruence of $\varrho$}. It is clear that $\varrho \subseteq \varrho ^{*}$.
Throughout the paper, the factor semigroup $S/\theta^{*}_S$ is denoted by $S/\theta^{*}$.

\medskip

In Section 2, we prove the following theorem.

\begin{theorem}\label{thm1}  Let $\beta$ be a congruence on a semigroup $S$ with $\theta _S\subseteq \beta$ and $P$ be a choice function on the collection of all $\theta_S$-classes of $S$. Then $\theta _{{\cal M}(S;S/\theta;P)}$ equals the kernel of ${\cal F}_{\theta _S, \beta}$ if and only if $\beta=\theta^{*}_S$. Thus ${\cal M}(S;S/\theta;P)/\theta \cong {\cal M}(S/\theta ^{*};S/\theta;P')$, where $P'$ is the mapping of $S/\theta$ onto $S/\theta^{*}$ defined as follows: $P'([s]_{\theta_S})=[s]_{\theta^{*}_S}$ for every $[s]_{\theta_S}\in S/\theta$.
\end{theorem}

\medskip

Let $\mathcal C$ be a class of semigroups. A congruence $\varrho$ on a semigroup $S$ is called a $\mathcal C$-congruence if the factor semigroup $S/\varrho$ belongs to $\mathcal C$. According to this definition, a congruence $\varrho$ on a semigroup $S$ is called a \emph{left cancellative congruence} if the factor semigroup $S/\varrho$ is left cancellative. In other words, $\varrho$ is a left cancellative congruence on a semigroup $S$ if $(xa,xb)\in \varrho$ implies $(a, b)\in \varrho$ for every $x, a, b\in S$.

In \cite[Theorem 1]{Cohn:sg-1}, P.M. Cohn gave necessary and sufficient conditions for a semigroup to be embedded in a left simple semigroup. The conditions differ essentially according to whether or not the semigroup contains an idempotent element. He proved that a semigroup $S$ is embedded in an idempotent-free left simple semigroup if and only if $S$ is idempotent-free and satisfies the condition: for all $a, b, x, y\in S$, $xa=xb$ implies $ya=yb$.
Using the terminology of \cite{Nagy:sg-8}, a semigroup $S$ satisfying this last condition is called a \emph{left equalizer simple semigroup}. In other words, a semigroup $S$ is left equalizer simple if, for arbitrary elements $a, b\in S$, the assumption that $xa=xb$ is satisfied for some $x\in S$ implies that $ya=yb$ is satisfied for all $y\in S$. By \cite[Theorem 2.1]{Nagy:sg-8}, a semigroup $S$ is left equalizer simple if and only if $\theta _S$ is a left cancellative congruence. Thus a semigroup $S$ is embedded in an idempotent-free left simple semigroup if and only if $S$ is idempotent-free and $\theta _S$ is a left cancellative congruence.

\medskip

In Section 3, we prove the following theorem.

\begin{theorem}\label{thm12} Let $S$ be a semigroup and $P$ be a choice function on the collection of all $\theta _S$-classes of $S$. Then the Rees matrix semigroup ${\mathcal M}(S;S/\theta ; P)$ can be embedded in an idempotent-free left simple semigroup if and only if $S$ is idempotent-free and the right colon congruence $\theta^{*}_S$ of $\theta_S$ is left cancellative.
\end{theorem}

\medskip

A semigroup satisfying the identity $axy=ayx$ is called a \emph{right commutative semigroup}.
As a corollary of Theorem~\ref{thm12}, we prove that if $S$ is an idempotent-free right commutative right cancellative semigroup and $P$ is a choice function on the collection of all $\theta_S$-classes of $S$, then the Rees matrix semigroup ${\mathcal M}(S;S/\theta;P)$ is embedded in an idempotent-free left simple semigroup (Corollary~\ref{cor12}). In addition, we show how to construct right commutative right cancellative semigroups (Theorem~\ref{thm122}).

\medskip

By the result of Theorem~\ref{thm1}, a finite sequence $S,\ S/\theta ,\ S/\theta ^{*}$ of semigroups has the property: there are mappings
\[S\stackrel{P}{\longleftarrow}S/\theta\stackrel{P'}{\longrightarrow}S/\theta ^{*}\]
such that
${\cal M}(S; S/\theta ; P)/\theta \cong {\cal M}(S/\theta ^{*}; S/\theta; P')$.
On the base of this result, we define a property of sequences of semigroups.

\begin{definition}
A triple $A, B, C$ of semigroups will said to be right regular, if there are mappings
\[A\stackrel{P}{\longleftarrow}B\stackrel{P'}{\longrightarrow}C\]
such that
${\cal M}(A; B; P)/\theta \cong {\cal M}(C; B; P')$.
An infinite sequence
\[S_0, S_1, \dots , S_{n-1}, S_n, S_{n+1}, \dots \]
of semigroups will said to be right regular if, for every positive integer $n$, the triple $S_{n-1}, S_n, S_{n+1}$ is right regular, i.e., there are mappings
\[S_{n-1}\stackrel{P_{n, n-1}}{\longleftarrow}S_n\stackrel{P_{n, n+1}}{\longrightarrow}S_{n+1}\]
such that
${\cal M}(S_{n-1}; S_{n}; P_{n, n-1})/\theta \cong {\cal M}(S_{n+1}; S_{n}; P_{n, n+1})$.
\end{definition}

\medskip
In Section 4, we show that every congruence on a semigroup defines a right regular sequence of semigroups.
For a congruence $\varrho$ on a semigroup $S$, we consider the sequence
$\varrho ^{(0)}\subseteq \varrho ^{(1)}\subseteq \cdots \subseteq \varrho ^{(n)}\subseteq \cdots$
of congruences, where $\varrho ^{(0)}=\varrho$ and, for every positive integer $n$, $\varrho ^{(n)}$ is the right colon congruence of $\varrho ^{(n-1)}$.
It is shown in \cite[Lemma 5]{Nagyleft:sg-6} that $\bigcup _{n=0}^{\infty}\varrho ^{(n)}\subseteq lrc(\varrho )$ for an arbitrary congruence $\varrho$ on an arbitrary semigroup $S$. The equation $\bigcup _{n=0}^{\infty}\varrho ^{(n)}=lrc(\varrho )$ is satisfied if and only if the congruence $\bigcup _{n=0}^{\infty}\varrho ^{(n)}$ is left reductive.
In Section 3, we show that sequences $\varrho ^{(0)}\subseteq \varrho ^{(1)}\subseteq \cdots \subseteq \varrho ^{(n)}\subseteq \cdots$ of congruences on semigroups also play an important role in our present investigation.
We prove the following theorem.

\begin{theorem}\label{thm2} For an arbitrary congruence $\varrho$ on an arbitrary semigroup $S$, the sequence
$S/\varrho ^{(0)},\ S/\varrho ^{(1)}, \dots ,\ S/\varrho ^{(n)}, \dots $
of factor semigroups is right regular, where $\varrho ^{(0)}=\varrho$ and $\varrho ^{(n)}$ ($n\geq 1$) denotes the right colon congruence of $\varrho ^{(n-1)}$.
\end{theorem}

For notations and notions not defined but used in this paper, we refer to books \cite{Clifford1:sg-1}, \cite{Howie:sg-11}, and \cite{Nagybook:sg-5}.

\section{The proof of Theorem~\ref{thm1}}

First we prove a lemma, which will be used in the proof of Theorem~\ref{thm1}.

\begin{lemma}\label{lem1} For arbitrary congruences $\alpha$ and $\beta$ on a semigroup $S$ with $\alpha \subseteq \beta$, and arbitrary choice function $P$ on the collection of all $\alpha$-classes of $S$, the mapping ${\cal F}_{\alpha , \beta}:(a, [s]_{\alpha})\mapsto ([a]_{\beta}, [s]_{\alpha})$ is a surjective homomorphism of the Rees matrix semigroup ${\cal M}(S; S/\alpha; P)$ onto the Rees matrix semigroup ${\cal M}(S/\beta ;S/\alpha ;P')$, where $P'$ is the mapping of $S/\alpha$ onto $S/\beta$ defined by $P'([s]_{\alpha})=[s]_{\beta}$ for every $[s]_{\alpha}\in S/\alpha$.
\end{lemma}

\noindent
{\bf Proof}. It is clear that ${\cal F}_{\alpha ,\beta}$ is surjective.
We show that ${\cal F}_{\alpha , \beta}$ is a homomorphism. Since $P([s]_{\alpha})\in [s]_{\alpha}$ and $\alpha \subseteq \beta$, we have
$[P([s]_{\alpha})]_{\beta}=[s]_{\beta}$ for every $s\in S$. Let
$(a, [s]_{\alpha})$ and $(b, [t]_{\alpha})$ be arbitrary elements of ${\cal M}(S; S/\alpha ; P)$. Then
\[{\cal F}_{\alpha , \beta}((a, [s]_{\alpha})(b, [t]_{\alpha}))={\cal F}_{\alpha, \beta}((aP([s]_{\alpha})b, [t]_{\alpha}))=\]
\[=([aP([s]_{\alpha})b]_{\beta}, [t]_{\alpha}))=([a]_{\beta}[P([s]_{\alpha})]_{\beta}[b]_{\beta}, [t]_{\alpha})=\]
\[=([a]_{\beta}[s]_{\beta}[b]_{\beta}, [t]_{\alpha})=([a]_{\beta}P'([s]_{\alpha})[b]_{\beta}, [t]_{\alpha})=\]
\[=([a]_{\beta}, [s]_{\alpha})([b]_{\beta}, [t]_{\alpha})={\cal F}_{\alpha , \beta}((a, [s]_{\alpha})){\cal F}_{\alpha , \beta}((b, [t]_{\alpha})).\]
Thus ${\cal F}_{\alpha , \beta}$ is a homomorphism.\hfill\openbox

\medskip

\noindent
{\bf The proof of Theorem~\ref{thm1}}:

\medskip
\noindent
We show that if $\beta$ is a congruence on a semigroup $S$ with $\theta _S\subseteq \beta$, then $\theta _{{\cal M}(S;S/\theta;P)}$ equals the kernel of ${\cal F}_{\theta _S, \beta}$ if and only if $\beta=\theta^{*}_S$.
Recall that $P$ is an arbitrary choice function on the collection of all $\theta _S$-classes of $S$ and ${\cal F}_{\theta _S, \beta}$ is the surjective homomorphism of the Rees matrix semigroup ${\cal M}(S; S/\theta; P)$ onto the Rees matrix semigroup ${\cal M}(S/\beta ;S/\theta ;P')$ defined as follows: ${\cal F}_{\theta _S, \beta}:(a, [s]_{\theta _S})\mapsto ([a]_{\beta}, [s]_{\theta _S})$ for every $a, s\in S$.

In the first part of the proof we show that $\theta _{{\cal M}(S;S/\theta;P)}$ equals the kernel of ${\cal F}_{\theta _S, \theta^{*}_S}$.
For elements $(a, [s]_{\theta _S})$ and $(b, [t]_{\theta _S})$ of  ${\cal M}(S; S/\theta ; P)$,
\[((a, [s]_{\theta _S}), (b, [t]_{\theta _S}))\in ker{\cal F}_{\theta _S, \theta ^{*}_S}\]
if and only if \[([a]_{\theta ^{*} _S}, [s]_{\theta _S}))=([b]_{\theta ^{*}_S}, [t]_{\theta _S}),\] i.e.,
\begin{equation}\label{111}
[s]_{\theta _S}=[t]_{\theta _S}
\end{equation} and
\begin{equation}\label{222}
(\forall x, y \in S)\ xy a=xy b.
\end{equation}
We show that (\ref{111}) and (\ref{222}) together are equivalent to the condition that
\[((a, [s]_{\theta _S}), (b, [t]_{\theta _S}))\in \theta _{{\cal M}(S;S/\theta :P)}.\]
Assume (\ref{111}) and (\ref{222}).
Then, for every $(x, [y]_{\theta _S})\in {\cal M}(S; S/\theta ; P)$, we have
\[(x, [y]_{\theta _S})(a, [s]_{\theta _S})=(xP([y]_{\theta _S})a, [s]_{\theta _S})=(xya, [s]_{\theta _S})= \]
\[=(xyb, [s]_{\theta _S})=(xP([y]_{\theta _S})b, [s]_{\theta _S})=(x, [y]_{\theta _S})(b, [t]_{\theta _S}),\] and hence
\[((a, [s]_{\theta _S}), (b, [t]_{\theta _S}))\in \theta _{{\cal M}(S;S/\theta :P)}.\]
Conversely, assume that
\[((a, [s]_{\theta _S}), (b, [t]_{\theta _S}))\in \theta _{{\cal M}(S;S/\theta :P)}.\]
Then, for every $x, y \in S$,
\[(x, [y]_{\theta _S})(a, [s]_{\theta _S})=(x, [y]_{\theta _S})(b, [t]_{\theta _S}),\] and hence
\[(xya, [s]_{\theta _S})=(xP([y]_{\theta _S})a, [s]_{\theta _S})=(x, [y]_{\theta _S})(a, [s]_{\theta _S})=\]
\[=(x, [y]_{\theta _S})(b, [t]_{\theta _S})=(xP([y]_{\theta _S})b, [t]_{\theta _S})=(xyb, [t]_{\theta _S}),\] which implies
(\ref{111}) and (\ref{222}).
Consequently \[ker{\cal F}_{\theta _S, \theta ^{*}_S}=\theta _{{\cal M}(S; S/\theta; P)}.\]

In the second part of the proof, we prove that if $ker{\cal F}_{\theta _S, \beta}=\theta _{{\cal M}(S; S/\theta; P)}$ for some congruence $\beta$ on $S$ with $\theta _S\subseteq \beta$, then $\beta =\theta ^{*}_S$. Assume \[ker{\cal F}_{\theta _S, \beta}=\theta _{{\cal M}(S; S/\theta; P)}\] for a congruence $\beta$ on $S$ with $\theta _S\subseteq \beta$.
First we show that $\beta \subseteq \theta ^{*}_S$. Assume $(a, b)\in \beta$, i.e., $[a]_{\beta}=[b]_{\beta}$, which is equivalent to the condition that, for every $[s]_{\theta _S}\in S/\theta$,
\[{\cal F}_{\theta _S, \beta}((a, [s]_{\theta _S}))={\cal F}_{\theta _S, \beta}((b, [s]_{\theta _S})),\] i.e., \[((a, [s]_{\theta _S}),(b, [s]_{\theta _S}))\in ker {\cal F}_{\theta _S, \beta}=\theta _{{\cal M}(S; S/\theta; P)}.\]
This last condition is equivalent to the condition that, for every $x, y \in S$, \[(x, [y ]_{\theta _S})(a, [s]_{\theta _S})=(x, [y ]_{\theta _S})(b, [s]_{\theta _S}),\] i.e.,
\[(xya, [s]_{\theta _S})=(xyb, [s]_{\theta _S}),\] which is equivalent to the condition that, for every $x, y \in S$, \[xy a=xy b,\] i.e.,
\[(a, b)\in \theta ^{*}_S.\] Hence \[\beta \subseteq \theta ^{*}_S.\]
Next we show that $\theta ^{*}_S\subseteq \beta$. Assume $(a, b)\in \theta ^{*}_S$. Then $xya=xyb$ for every $x, y\in S$. Let $(x, [y]_{\theta _S})$ be an arbitrary element of ${\cal M}(S; S/\theta; P)$, and let $[s]_{\theta _S}$ be an arbitrary element of $S/\theta$. Then
 \[(x, [y]_{\theta _S})(a, [s]_{\theta _S})=(xya, [s]_{\theta _S})=(xyb, [s]_{\theta _S})=(x, [y]_{\theta _S})(b, [s]_{\theta _S}),\]
i.e., \[((a, [s]_{\theta _S}),(b, [s]_{\theta _S}))\in \theta _{{\cal M}(S;S/\theta ;P)}=ker {\cal F}_{\theta _S, \beta}.\]
Then $[a]_{\beta}=[b]_{\beta}$, i.e., $(a, b)\in \beta$. Hence \[\theta ^{*}_S\subseteq \beta.\] Consequently $\beta =\theta ^{*}_S$.

The isomorphism ${\cal M}(S;S/\theta;P)/\theta \cong {\cal M}(S/\theta ^{*}_S;S/\theta;P')$ is a consequence of the equation
$ker {\cal F}_{\theta _S, \theta ^{*}_S}=\theta _{{\cal M}(S;S/\theta;P)}$.
\hfill\openbox

\section{The proof of Theorem~\ref{thm12}}

We show that if $S$ is a semigroup and $P$ is a choice function on the collection of all $\theta _S$-classes of $S$, then the Rees matrix semigroup ${\mathcal M}(S;S/\theta ; P)$ can be embedded in an idempotent-free left simple semigroup if and only if $S$ is idempotent-free and the right colon congruence $\theta^{*}_S$ of $\theta_S$ is left cancellative.

By \cite[Corollary 2.4]{Nagy:sg-21},
a semigroup $S$ is idempotent-free if and only if the factor semigroup $S/\theta$ is idempotent-free. This fact will be used throughout the proof.

We assume first that $S$ is an idempotent-free semigroup such that the right colon congruence $\theta ^{*}_S$ is left cancellative. Then the factor semigroup $S/\theta^{*}$ is left cancellative. Since $S$ is idempotent-free and $S/\theta^{*}\cong (S/\theta)/\theta$, the factor semigroup $S/\theta^{*}$ is idempotent-free. Thus $S/\theta^{*}$ is left cancellative and idempotent-free. We show that the Rees matrix semigroup
${\mathcal M}(S/\theta ^{*}; S/\theta; P')$ is also left cancellative and idempotent-free, where $P': S/\theta \mapsto S/\theta ^{*}$ is defined by $P'([s]_{\theta})=[s]_{\theta^{*}}$. Let $([x]_{\theta ^{*}_S}, [y]_{\theta_S})$, $([a]_{\theta ^{*}_S}, [s]_{\theta_S})$, $([b]_{\theta ^{*}_S}, [t]_{\theta_S})$ be arbitrary elements of ${\mathcal M}(S/\theta ^{*}; S/\theta; P')$. Assume
\[([x]_{\theta ^{*}_S}, [y]_{\theta_S})([a]_{\theta ^{*}_S}, [s]_{\theta_S})=([x]_{\theta ^{*}_S}, [y]_{\theta_S})([b]_{\theta ^{*}_S}, [t]_{\theta_S}).\] Then
\[([xya]_{\theta ^{*}_S}, [s]_{\theta_S})=([x]_{\theta^{*}_S}[y]_{\theta^{*}_S}[a]_{\theta ^{*}_S}, [s]_{\theta_S})=
([x]_{\theta^{*}_S}P'([y]_{\theta_S})[a]_{\theta ^{*}_S}, [s]_{\theta_S})=\]
\[=([x]_{\theta ^{*}_S}, [y]_{\theta_S})([a]_{\theta ^{*}_S}, [s]_{\theta_S})=([x]_{\theta ^{*}_S}, [y]_{\theta_S})([b]_{\theta ^{*}_S}, [t]_{\theta_S})=\]
\[=([x]_{\theta^{*}_S}P'([y]_{\theta_S})[b]_{\theta ^{*}_S}, [t]_{\theta_S})
=([x]_{\theta^{*}_S}[y]_{\theta^{*}_S}[b]_{\theta ^{*}_S}, [t]_{\theta_S})=[xyb]_{\theta ^{*}_S}, [t]_{\theta_S}),\]
and hence
\[(xya, xyb)\in \theta^{*}_S\quad \hbox{and}\quad [s]_{\theta_S}=[t]_{\theta_S}.\] Since $\theta^{*}_S$ is a left cancellative congruence on $S$, we have $(a, b)\in \theta^{*}_S$, and consequently  \[([a]_{\theta ^{*}_S}, [s]_{\theta_S})=([b]_{\theta ^{*}_S}, [t]_{\theta_S}).\] Thus the Rees matrix semigroup ${\mathcal M}(S/\theta ^{*}; S/\theta; P')$ is left cancellative. Suppose, indirectly, that the Rees matrix semigroup
${\mathcal M}(S/\theta ^{*}; S/\theta; P')$ contains an idempotent element $([e]_{\theta ^{*}_S}, [f]_{\theta_S})$. Then
\[([e]_{\theta ^{*}_S}[f]_{\theta ^{*}_S}[e]_{\theta ^{*}_S}, [f]_{\theta_S})=([e]_{\theta ^{*}_S}, [f]_{\theta_S}),\] and hence $[e]_{\theta ^{*}_S}$ is a regular element of the semigroup $S/\theta^{*}$. Thus $S/\theta^{*}$ contains an idempotent element. This is a contradiction. Consequently the Rees matrix semigroup ${\mathcal M}(S/\theta ^{*}; S/\theta; P')$ is an idempotent-free left cancellative semigroup. By Theorem~\ref{thm1},
${\cal M}(S;S/\theta;P)/\theta \cong {\cal M}(S/\theta ^{*};S/\theta;P')$. Thus $\theta_{{\mathcal M}(S;S/\theta ;P)}$ is a left cancellative congruence, and hence ${\cal M}(S;S/\theta;P)$ is left equalizer simple by \cite[Theorem 2.1]{Nagy:sg-8}. Since ${\cal M}(S/\theta ^{*};S/\theta;P')$ is idempotent-free, ${\cal M}(S;S/\theta;P)/\theta \cong {\cal M}(S/\theta ^{*};S/\theta;P')$ also implies that ${\cal M}(S;S/\theta;P)$ is idempotent-free. Thus ${\cal M}(S;S/\theta;P)$ is an idempotent-free left equalizer simple semigroup. By \cite[Theorem 1]{Cohn:sg-1}, ${\cal M}(S;S/\theta;P)$ can be embedded in an idempotent-free left simple semigroup.

Conversely, assume that the Rees matrix semigroup ${\cal M}(S;S/\theta;P)$ can be embedded in an idempotent-free left simple semigroup. Then
${\cal M}(S;S/\theta;P)$ is an idempotent-free left equalizer simple semigroup by \cite[Theorem 1]{Cohn:sg-1}. By Theorem~\ref{thm1},
${\cal M}(S/\theta ^{*};S/\theta;P')\cong {\cal M}(S;S/\theta;P)/\theta$. Then, applying also \cite[Theorem 2.1]{Nagy:sg-8}, we get that ${\cal M}(S/\theta ^{*};S/\theta;P')$ is an idempotent-free left cancellative semigroup. If $e$ is an idempotent element of $S$, then
\[([e]_{\theta^{*}_S}, [e]_{\theta _S})([e]_{\theta^{*}_S}, [e]_{\theta _S})=([e]_{\theta^{*}_S}, [e]_{\theta _S}),\] i.e., $([e]_{\theta^{*}_S}, [e]_{\theta _S})$ is an idempotent element of ${\cal M}(S/\theta ^{*};S/\theta;P')$. Thus $S$ is idempotent-free, because ${\cal M}(S/\theta ^{*};S/\theta;P')$ is idempotent-free. It remains to show that $\theta ^{*}_S$ is a left cancellative congruence on $S$. Assume $(xa,xb)\in \theta ^{*}_S$ for elements $x, a, b\in S$. Then, for arbitrary $y, s\in S$, we have
\[([y]_{\theta ^{*}_S}, [x]_{\theta_S})([a]_{\theta ^{*}_S}, [s]_{\theta_S})=
([y]_{\theta^{*}_S}P'([x]_{\theta_S})[a]_{\theta^{*}_S}, [s]_{\theta_S})=([yxa]_{\theta ^{*}_S}, [s]_{\theta_S})=\]
\[=([yxb]_{\theta ^{*}_S}, [s]_{\theta_S})=([y]_{\theta^{*}_S}P'([x]_{\theta_S})[b]_{\theta^{*}_S}, [s]_{\theta_S})=
([y]_{\theta ^{*}_S}, [x]_{\theta_S})([b]_{\theta ^{*}_S}, [s]_{\theta_S}),\] from which we get
\[([a]_{\theta ^{*}_S}, [s]_{\theta_S})=([b]_{\theta ^{*}_S}, [s]_{\theta_S}),\]
because ${\cal M}(S/\theta ^{*};S/\theta;P')$ is left cancellative. Thus $(a,b)\in \theta ^{*}_S$. Consequently $\theta ^{*}_S$ is a left cancellative congruence.\hfill\openbox

\begin{corollary}\label{cor12} Let $S$ be an idempotent-free right commutative right cancellative semigroup and $P$ be a choice function on the collection of all $\theta_S$-classes of $S$. Then the Rees matrix semigroup ${\mathcal M}(S;S/\theta;P)$ is embedded in an idempotent-free left simple semigroup.
\end{corollary}

\noindent
{\bf Proof}. By Theorem~\ref{thm12}, it is sufficient to show that $\theta^{*}_S$ is a left cancellative congruence on $S$. Assume $(xa, xb)\in \theta^{*}_S$ for elements $x, a, b\in S$. Then, for every $s, t\in S$, $stxa=stxb$, and hence $stax=stbx$ by right commutativity of $S$. Since $S$ is right cancellative, $sta=stb$ (for every $s, t\in S$), and hence $(a, b)\in \theta^{*}_S$. Consequently $\theta^{*}_S$ is a left cancellative congruence.\hfill\openbox

\medskip

\begin{corollary}\label{cor22} The following conditions on a right commutative right cancellative semigroup $S$ are equivalent.
\begin{itemize}
\item[(1)] $S$ is embedded in an idempotent-free left simple semigroup.
\item[(2)] For any choice function $P$ on the collection of all $\theta_S$-classes of $S$, the Rees matrix semigroup ${\mathcal M}(S;S/\theta ;P)$ is embedded in an idempotent-free left simple semigroup.
\end{itemize}
\end{corollary}

\noindent
{\bf Proof}.
$(1)$ implies $(2)$: If $S$ is embedded in an idempotent-free left simple semigroup, then $S$ is idempotent-free, and hence the Rees matrix semigroup ${\mathcal M}(S;S/\theta ;P)$ is embedded in an idempotent-free left simple semigroup by Corollary~\ref{cor12}.

$(2)$ implies $(1)$: Assume that the Rees matrix semigroup ${\mathcal M}(S;S/\theta ;P)$ is embedded in an idempotent-free left simple semigroup. Then, by Theorem~\ref{thm12}, $S$ is idempotent-free and $\theta^{*}_S$ is a left cancellative congruence. Assume $(a, b)\in \theta^{*}_S$ for elements
$a, b\in S$. Then, for all $x, y\in S$, $xya=xyb$. Since $S$ is right commutative, $xay=xby$. Since $S$ is right cancellative, we get $xa=xb$ (for all $x\in S$). Hence $(a, b)\in \theta _S$. Thus $\theta^{*}_S\subseteq \theta_S$, and consequently $\theta _S=\theta^{*}_S$, because
$\theta_S\subseteq \theta^{*}_S$ is satisfied by the definition of the right colon congruence. Hence $\theta_S$ is a left cancellative congruence. Thus, by \cite[Theorem 2.1]{Nagy:sg-8}, $S$ is left equalizer simple. Consequently $S$ is embedded in an idempotent-free left simple semigroup by \cite[Theorem 1]{Cohn:sg-1}.\hfill\openbox

\medskip

In the next we show how to construct (idempotent-free) right commutative right cancellative semigroups. We use the following construction.

\begin{construction}\rm (\cite[Construction 1]{Nagy:sg-8})\label{constr}
Let $T$ be a left cancellative semigroup. For each $t\in T$, associate a nonempty set $S_t$ to $t$ such that $S_t\cap S_r=\emptyset$ for every
$t, r\in T$ with $t\neq r$.

For arbitrary couple $(t, r)\in T\times T$ with $r\in tT$, let $\varphi_{t, r}$ be a mapping of $S_t$ into $S_r$ acting on the right.
For all $t\in T$, $r\in tT$, $q\in rT\subseteq tT$, assume \[\varphi_{t,r}\varphi_{r,q}=\varphi_{t, q}.\]

On the set $S=\cup _{t\in T}S_t$ define an operation $*$ as follows: for arbitrary $a\in S_t$ and $b\in S_x$, let \[a* b=a\varphi_{t, tx}.\]

If $a\in S_t$, $b\in S_x$, $c\in S_y$ are arbitrary elements then
\[a* (b* c)=a* b\varphi_{x, xy}=a\varphi_{t, t(xy)}=\]
\[=a\varphi_{t, tx}\varphi_{tx, t(xy)}=a\varphi_{t, tx}* c=(a* b)* c.\] Thus
$(S; * )$ is a semigroup.
\end{construction}

\medskip
The semigroup $(S; * )$ defined in Construction~\ref{constr} will said to be a \emph{right regular extension} of the left cancellative semigroup $T$.

A right regular extension of a left cancellative semigroup $T$ is is said to be \emph{injective} if the mappings $\varphi_{t, r}$ are injective for every $t\in T$ and $r\in tT$.

\begin{theorem}\label{thm122} A semigroup $S$ is right commutative and right cancellative if and only if it is an injective right regular extension of a commutative cancellative semigroup $T$. $S$ is idempotent-free if and only if $T$ is idempotent-free.
\end{theorem}

\noindent
{\bf Proof}. Let $S$ be a right commutative right cancellative semigroup. Assume $xa=xb$ for elements $x, a, b\in S$. Then \[uax=uxa=uxb=ubx\] is satisfied for all $u\in S$. Since $S$ is right cancellative, we get $ua=ub$ for all $u\in S$. Hence $S$ is left equalizer simple. Then, by \cite[Theorem 2.2]{Nagy:sg-8}, $S$ is a right regular extension of a left cancellative semigroup $T$. Since $S$ is right commutative, $T$ is commutative. Thus $S$ is a right regular extension of the commutative cancellative semigroup $T$. We show that the mappings $\varphi_{t, r}$ are injective for every $t\in T$ and $r\in tT$. Let $t\in T$, $r\in tT$, and $a,b\in S_t$ be arbitrary elements. Assume $a\varphi_{t, r}=b\varphi_{t, r}$. Let $x\in T$ be the (only) element of $T$ with $r=tx$. Denoting the multiplication on $S$ by $*$, for an arbitrary $c\in S_x$, \[a*c=a\varphi_{t, r}=b\varphi_{t, r}=b*c,\] from which we get $a=b$, because $S$ is right cancellative. Thus $\varphi_{t, r}$ is injective. Consequently $S$ is an injective right regular extension of the commutative cancellative semigroup $T$.

Conversely, assume that a semigroup $S$ is an injective right regular extension of a commutative cancellative semigroup $T$. It is clear that $S$ is right commutative. To show that $S$ is right cancellative, assume $a*c=b*c$ for elements $a\in S_x$, $b\in S_y$, $c\in S_t$ ($x, y, t\in T$). Then $a\varphi_{x,xt}=b\varphi_{y,yt}$, and hence $xt=yt$. Since $T$ is cancellative, $x=y$. Thus $a, b\in S_x$ and $a\varphi_{x,xt}=b\varphi_{x, xt}$. Since $\varphi_{x, xt}$ is injective, we get $a=b$. Consequently $S$ is right cancellative.

Since $T\cong S/\theta$, the proof of Theorem~\ref{thm12} implies that $S$ is idempotent-free if and only if $T$ is idempotent-free.\hfill\openbox

\section{The proof of Theorem~\ref{thm2}}

For brevity, the factor semigroup $S/\theta ^{(k)}_S$ on a semigroup $S$ will be denoted by $S/\theta ^{(k)}$ for every nonnegative integer $k$.
First we prove a lemma which will be used in the proof of Theorem~\ref{thm2}.

\begin{lemma}\label{kieg} Let $S$ be a semigroup. Then, for every non-negative integers $n$ and $k$, $(S/\varrho ^{(n)})/\theta ^{(k)}\cong S/\varrho ^{(n+k+1)}$.
\end{lemma}

\noindent
{\bf Proof}. Let $S$ be an arbitrary semigroup and $n$, $k$ be arbitrary non-negative integers.
Let $(\cdot )\varphi _1$ denote the canonical homomorphism of $S$ onto $S/\varrho ^{(n)}$, and let $(\cdot )\varphi _2$ be the canonical homomorphism of $(S/\varrho ^{(n)})/\theta ^{(k)}$. We show that the kernel of $\varphi _1\circ \varphi _2$ is $\varrho ^{(n+k+1)}$, and so the assertion of the lemma follows from the homomorphism theorem.

Let $a, b\in S$ be arbitrary elements. Then
\[(a)(\varphi _1\circ \varphi _2)=(b)(\varphi _1\circ \varphi _2)\]
if and only if
\[([a]_{\varrho ^{(n)}}, [b]_{\varrho ^{(n)}})\in \theta ^{(k)}_{S/\varrho ^{(n)}},\] which is equivalent to the condition that, for every $x_i\in S$ ($i=1, \dots , k+1$),
\[[x_1]_{\varrho ^{(n)}}\cdots [x_{k+1}]_{\varrho ^{(n)}}[a]_{\varrho ^{(n)}}=[x_1]_{\varrho ^{(n)}}\cdots [x_{k+1}]_{\varrho ^{(n)}}[b]_{\varrho ^{(n)}},\]
i.e.,
\[(x_1\cdots x_{k+1}a, x_1\cdots x_{k+1}b)\in \varrho ^{(n)},\]
which is equivalent to the condition that
\[(a, b)\in \varrho ^{n+k+1}.\]
Consequently, the kernel of $\varphi _1\circ \varphi _2$ equals $\varrho ^{(n+k+1)}$.\hfill\openbox

\bigskip

\noindent
{\bf The proof of Theorem~\ref{thm2}}:

\medskip

Let $\varrho$ be an arbitrary congruence on a semigroup $S$. We show that the sequence \[S/\varrho ^{(0)},\ S/\varrho ^{(1)}, \dots ,\ S/\varrho ^{(n-1)},\ S/\varrho ^{(n)},\ S/\varrho ^{(n+1)}, \dots \]
of factor semigroups is right regular. Let $n$ be an arbitrary positive integer.
Applying Lemma~\ref{kieg}, we have
\[(S/\varrho ^{(n-1)})/\theta \cong S/\varrho ^{(n)}\quad \hbox{and}\quad (S/\varrho ^{(n-1)})/\theta ^{(1)}\cong  S/\varrho ^{(n+1)}.\] Then, by Theorem~\ref{thm1}, the triple \[S/\varrho ^{(n-1)}, S/\varrho ^{(n)}, S/\varrho ^{(n+1)}\] of factor semigroups is right regular. Consequently, the sequence
\[ S/\varrho ^{(0)}, \dots, S/\varrho ^{(n-1)}, S/\varrho ^{(n)}, S/\varrho ^{(n+1)}, \dots \]
is right regular. \hfill\openbox

\begin{remark}\rm
By Lemma~\ref{kieg}, $(S/\theta)/\theta\cong S/\theta ^{*}$ is satisfied for an arbitrary semigroup $S$. Thus, if $A, B, C$ are semigroups such that
$A/\theta \cong B$ and $B/\theta \cong C$, then the triple $A, B, C$ is right regular by Theorem~\ref{thm1}.
In the next, we give a right regular triple $A, B, C$ of semigroups, where non of the isomorphisms $A/\theta \cong B$ and $B/\theta \cong C$ is fulfilled.

Let $A$ be a left zero semigroup, i.e., a semigroup which satisfies the identity $ab=a$. It is easy to see that $\theta _A$ is the universal relation on $A$. Let $B$ be a right zero semigroup, i.e., a semigroup which satisfies the identity $ab=b$. It is clear that $\theta _B$ is the identity relation on $B$. Assume $|B|\geq 2$. Then $A/\theta \not \cong B$. If $C=\{e\}$ is a one-element semigroup, then $B/\theta \not \cong C$. Thus non of the isomorphisms $A/\theta \cong B$ and $B/\theta \cong C$ is fulfilled. We show that the triple
$A, B, C$ is right regular.
Let $P: B\to A$ be an arbitrary mapping. For every elements
$(a_1, b_1)$ and $(a_2, b_2)$ of ${\cal M}(A; B; P)$, $(a_1, b_1)(a_2, b_2)=(a_1P(b_1)a_2, b_2)=(a_1, b_2)$.
Thus ${\cal M}(A; B; P)$ is isomorphic to the direct product of the left zero semigroup $A$ and the right zero semigroup $B$. It is easy to see that
${\cal M}(A; B; P)/\theta \cong B$.
Let $P'$ denote the (only possible) mapping of $B$ onto the one-element semigroup $C=\{ e\}$. In the semigroup ${\cal M}(C; B; P')$, we have
$(e, b_1)(e, b_2)=(eP'(b_1)e, b_2)=(e, b_2)$
for every $b_1, b_2\in B$ from which it follows that $b\mapsto (e, b)$ is an isomorphism of $B$ onto ${\cal M}(C; B; P')$, and so
$B\cong {\cal M}(C; B; P')$. Consequently
${\cal M}(A; B; P)/\theta \cong {\cal M}(C; B; P')$. Hence
the triple $A, B, C$ is right regular.
\end{remark}

\end{document}